\documentclass[a4paper, 10pt, twoside]{article}

\usepackage{exscale, fullpage}
\usepackage[centertags]{amsmath}
\usepackage{amssymb}
\usepackage{amsthm}
\usepackage{dsfont}
\usepackage{url}
\usepackage[hypertex]{hyperref}

\newtheorem{definition}{Definition}

\newtheorem{theorem}[definition]{Theorem}
\newtheorem{proposition}[definition]{Proposition}
\newtheorem{lemma}[definition]{Lemma}
\newtheorem{corollary}[definition]{Corollary}
\newtheorem{conjecture}[definition]{Conjecture}

\newtheorem{deflemma}[definition]{Definition-Lemma}

\newcommand{\nd}{\noindent}

\newcommand{\dC}{{\mathds C}}
\newcommand{\dQ}{{\mathds Q}}

\newcommand{\dZ}{{\mathds Z}}
\newcommand{\dP}{{\mathds P}}

\newcommand{\dH}{{\mathbb H}}
\newcommand{\dL}{{\mathbb L}}

\newcommand{\bR}{{\mathbb R}}

\newcommand{\cD}{\mathcal{D}}
\newcommand{\cE}{\mathcal{E}}

\newcommand{\cG}{\mathcal{G}}
\newcommand{\cH}{\mathcal{H}}
\newcommand{\cI}{\mathcal{I}}

\newcommand{\cK}{\mathcal{K}}
\newcommand{\cL}{\mathcal{L}}
\newcommand{\cM}{\mathcal{M}}
\newcommand{\cN}{\mathcal{N}}
\newcommand{\cO}{\mathcal{O}}

\newcommand{\cR}{\mathcal{R}}

\DeclareMathOperator{\Spec}{\textup{Spec}\,}

\DeclareMathOperator{\Sp}{\textup{Sp}_{\theta=0}}

\DeclareMathOperator{\diag}{\textup{diag}}
\DeclareMathOperator{\vol}{\textup{vol}}

\DeclareMathOperator{\Gl}{Gl}

\DeclareMathOperator{\Lie}{\textup{Lie}}

\DeclareMathOperator{\DR}{\mathit{DR}}
\newcommand{\Weyl}{{\dC[t]\langle\partial_t\rangle}}

\begin{document}

\title{Bernstein polynomials and spectral numbers \\ for linear free divisors}
\author{Christian Sevenheck}

\maketitle

\begin{abstract}
We discuss Bernstein polynomials of reductive linear free divisors.
We define suitable Brieskorn lattices for these non-isolated singularities,
and show the analogue of Malgrange's result relating the roots of
the Bernstein polynomial to the residue eigenvalues on the saturation
of these Brieskorn lattices.
\end{abstract}
\renewcommand{\thefootnote}{}
\footnote{
This research is partially supported by ANR grant ANR-08-BLAN-0317-01 (SEDIGA).
}

\section{Introduction}
\label{sec:Introduction}

In this note, we show that for reductive linear free divisors $D\subset \dC^n$,
which were studied in a number of recent papers (see \cite{BM}
\cite{GMNS}, \cite{dGMS} and \cite{GS}), the roots of the Bernstein
polynomial of a defining equation $h$ of $D$ can be recovered as a
certain set of eigenvalues of a residue endomorphism. More precisely,
for a generic linear form $f$ on $\dC^n$, one defines
a family of Gau\ss-Manin systems for $f$, seen as a function
on the fibres of $h$. This family has a specific (logarithmic) extension
over $D$, which gives the set of residue eigenvalues we are interested in.

This relation to the Bernstein polynomial
has a number of consequences: First, while the definition of
the Bernstein polynomial is rather simple, it is in general very hard
to calculate its roots in concrete examples. This is true even for linear free divisors,
though the differential operator occurring in Bernstein's functional
equation is more explicitly known than in the general case.
On the other hand, the calculation of the residue eigenvalues alluded
to above is, although not trivial, easier to carry out. We apply
our result to obtain, using the calculations from \cite[chapter 6]{dGMS},
Bernstein polynomials for discriminants in
representation spaces of the Dynkin quivers $A_n$, $D_n$ and $E_6$
as well as the so-called star quiver $\star_n$, also considered in loc.cit. (which is
not a Dynkin quiver for $n>3$). We also calculate the Bernstein polynomials
for two irreducible linear free divisors which are discriminants
of irreducible pre-homogenous vector spaces described in \cite{sk}.

Another motivation for this work comes from the fact that the residue
eigenvalues of the Gau\ss-Manin systems give information on their limit
behavior (resp., of the corresponding family of Brieskorn lattices), when approaching
the zero fibre of $h$, i.e., the divisor $D$. In particular,
in \cite{dGMS}, questions about the degeneration of \emph{Frobenius manifolds},
associated to the tame functions $f_{|D_t}$, where $D_t:=h^{-1}(t)$, were related to the asymptotic behavior of
a natural pairing defined on the Gau\ss-Manin system. In particular, the residue
eigenvalues of these Gau\ss-Manin systems then need to be symmetric around zero.
This was stated as a conjecture in loc.cit., and it follows from the relation
between these eigenvalues and the roots of the Bernstein polynomial that we prove here.

Finally, the family of Brieskorn lattices associated to $f_{|D_t}$ also has logarithmic
extension over the divisor $D$, constructed using logarithmic differential forms.
We define the fibre over $t=0$ to be the \emph{logarithmic Brieskorn lattice} of
$D$ (in fact, it does not depend on the choice of the linear form).
In contrast to the fibres at $t\neq 0$, this Brieskorn lattice
is regular singular at the origin, reflecting the local situation of the
pair $(f,h)$ at the origin in $\dC^n$. It turns out that
then our result can be rephrased to give the analogue of Malgrange's classical
result for isolated singularities: The roots of the Bernstein polynomial
are (up to a rescaling) the residue eigenvalues on the \emph{saturation} of
this logarithmic Brieskorn lattice.

\textbf{Acknowledgements:} I thank Claude Sabbah and Michel Granger for their
help during the preparation of this article and Mathias Schulze for
comments on a first version. I am particularly
grateful to Ignacio de Gregorio for all the discussions on linear free divisors
and related subjects we had over the last years and also for
having done a good part of the computations used in this paper.

\section{Linear Free Divisors and Gau\ss-Manin systems}
\label{sec:LFDandGM}

In this section, we first recall from \cite{dGMS} the construction of the family
of Gau\ss-Manin systems associated to a linear section of a linear
free divisor. We also give a more intrinsic definition of these
Gau\ss-Manin systems as a direct image of a map constructed from two
polynomials. Finally, we discuss the definition of the residue eigenvalues
relevant for the present work, as introduced in \cite{dGMS}.

We denote throughout this article by $V$ the complex vector space $\dC^n$.
\begin{deflemma}[\cite{KS1}, \cite{BM}, \cite{dGMS}]
\begin{enumerate}
\item
Let $D\subset V$ be a reduced hypersurface with defining equation $h$. Then $D$ is called a free divisor,
if the sheaf $\Theta_V(-\log\,D):=\{\vartheta\in\Theta_V\,|\,\vartheta(h)\subset(h)\}$ is
a free $\cO_V$-module. If moreover a basis $(\xi_i)$ of $\Theta_V(-\log\, D)$ exists
such that $\xi_i=\sum_{j=1}^n \xi_{ij}\partial_{x_j}$ where $\xi_{ij}$ are \textbf{linear} forms
on $V$, then $D$ is called linear free.
\item
Let $G$ be the identity component of the algebraic group
$G_D:=\left\{g\in\Gl(V)\,|\,g(D)\subset D\right\}$. Then $(V,G)$ is a
pre-homogenous vector space in the sense of Sato (see, e.g., \cite{sk}),
in particular, the complement $V\backslash D$ is an open orbit of $G$.
We call $D$ reductive if $G_D$ is so. A rational function $r\in\dC(V)$
is called a semi-invariant if there is a character $\chi_r:G\rightarrow \dC^*$ such
that $g(r)=\chi_r(g)\cdot r$ for all $g\in G$. Obviously, $h$ itself is a semi-invariant.
\item
$G$ acts on $V^*$ by the dual action, with dual discriminant $D^*\subset V^*$.
If $G_D$ is reductive, then $(V^*,D^*)$ is pre-homogenous.
We call a linear form $f\in V^*$ generic with respect to $h$ (or simply generic, if
no confusion is possible) if $f$ lies in the open orbit
$V^*\backslash D^*$ of the dual action.

There is a basis $(e_i)$ of $V$ with corresponding coordinates $(x_i)$ (called unitary) such
that $G$ appears as a subgroup of $U(n)$ in these coordinates.
Then $D^*=\{h^*=0\}$, where $h^*(y):=\overline{h(\overline{y})}$,
$(y_i)$ being the dual coordinates of $(x_i)$.
\end{enumerate}
\end{deflemma}

In the sequel, we always consider linear forms which are generic
with respect to $h$.

In order to study the behavior of the restriction of the linear function $f$ on the fibres
$D_t:=h^{-1}(t)$, $t \neq 0$, but also on $D$ itself, the following \emph{deformation algebra}
was introduced in \cite{dGMS}.
\begin{definition}
Let $D$ be linear free with defining equation $h$, seen as a morphism
$h:V\rightarrow T:=\Spec \dC[t]$.
\begin{enumerate}
\item
Let $E\in \Theta_V(-\log\,D)$ be the Euler field $E=\sum_{i=1}^n x_i\partial_{x_i}$.
Call
$$
\Theta_{V/T}(-\log\,D):=\{\vartheta\in\Theta_V(-\log\,D)\,|\,\vartheta(h)=0\}
$$
the module of relative logarithmic vector fields. $\Theta_{V/T}(-\log\,D)$ is $\cO_V$-free
of rank $n-1$, and we have a decomposition $\Theta_V(-\log\,D)=\cO_V E \oplus \Theta_{V/T}(-\log\,D)$.
\item
The ideal $J_h(f):=df(\Theta_{V/T}(-\log\,D)) \subset \cO_V$ is called the Jacobian ideal of the pair $(f,h)$.
The quotient $\cO_V/J_h(f)$ is the Jacobian algebra (or deformation algebra) of $(f,h)$.
\end{enumerate}
\end{definition}
Notice that $\Theta_{V/T}(-\log\,D)$ was called $\Theta_V(-\log\,h)$ in \cite{dGMS}. It was shown in loc.cit.,
section 3.2, that if $f$ is generic with respect to $h$, then $h_*\cO_V/J_h(f)$ is $\cO_T$-free of rank $n$,
and generated by $(f^i)$ for $i=0,\ldots,n-1$. Moreover,
it is interpreted as the relative tangent space $T^1_{\cR_h/\dC}(f)$ of the deformation theory of $f$ with respect to the
group $\cR_h$ of right-equivalences preserving all fibres of $h$.

Denote by $(\Omega^\bullet_{V/T}(\log\,D),d):=(\Omega^\bullet_V(\log\,D)/(\Omega^{\bullet-1}_V(\log\,D)\wedge h^*\Omega^1(\log\,\{0\}),d)$ the relative logarithmic de Rham complex of $h$ as studied, under the name
$\Omega^\bullet(\log\,h)$ in \cite[section 2.2]{dGMS}). This relative logarithmic complex is used in the definition
of the family of Gau\ss-Manin-systems resp. Brieskorn lattices in loc.cit.
\begin{deflemma}[{\cite[section 4]{dGMS}}]
Let $h$ and $f$ as above. Define
\begin{equation}
\label{eq:GM-Brieskorn}
\begin{array}{rclcrcl}
G(\log\,D) & := & \frac{H^0(V,\Omega^{n-1}_{V/T}(\log\,D)[\theta,\theta^{-1}])}{(\theta d-df\wedge)H^0(V,\Omega^{n-2}_{V/T}(\log \,D)[\theta,\theta^{-1}])}
&;&
G(*D) & := & \frac{H^0(V,\Omega^{n-1}_{V/T}(*D)[\theta,\theta^{-1}])}{(\theta d-df\wedge)H^0(V,\Omega^{n-2}_{V/T}(*D)[\theta,\theta^{-1}])}
\\ \\
G_0(\log\,D) & := & \frac{H^0(V,\Omega^{n-1}_{V/T}(\log \,D)[\theta])}{(\theta d-df\wedge)H^0(V,\Omega^{n-2}_{V/T}(\log \,D)[\theta])}
&;&
G_0(*D) & := & \frac{H^0(V,\Omega^{n-1}_{V/T}(*D)[\theta])}{(\theta d-df\wedge)H^0(V,\Omega^{n-2}_{V/T}(*D)[\theta])}.
\\ \\
\end{array}
\end{equation}
Then $G(*D)$ is $\dC[\theta,\theta^{-1},t,t^{-1}]$-free of rank $n$ and
$G(\log\,D)$ (resp. $G_0(*D), G_0(\log\, D)$) is a
$\dC[\theta,\theta^{-1},t]$- (resp. $\dC[\theta,t,t^{-1}]$-,
$\dC[\theta,t]$-) lattice inside $G(*D)$. These modules fit into the following diagram
$$
\begin{array}{ccc}
G(\log\,D) & \subset & G(*D) \\
\cup & &\cup \\
G_0(\log\, D) &\subset & G_0(*D).
\end{array}
$$
Define a connection
$$
\nabla:G_0(\log\,D)\longrightarrow G_0(\log\,D)\otimes \theta^{-1}\Omega^1_{\dC\times T}\left(\log(\{0\}\times T)\cup(\dC\times \{0\})\right)
$$
by putting, for a form $\omega\in H^0(V,\Omega^{n-1}_{V/T}(\log\,D))$,
\begin{equation}\label{eq:ConnectionOperator}
\begin{array}{rcl}
\nabla_{\partial_\theta}([\omega]) & := & \theta^{-2}[f\cdot\omega] \\ \\
\nabla_{\partial_t}([\omega]) & := & \frac1{nt}([\Lie_E(\omega)]-[\theta^{-1} f\cdot\omega])
\end{array}
\end{equation}
and extending by the Leibniz-rule (for $\nabla_{\partial_\theta}$) resp. $\theta$-linearly (for
$\nabla_{\partial_t}$). We denote by $\nabla$ the induced connection on
$G(\log\,D)$, $G_0(*D)$ and $G(*D)$.
\end{deflemma}

One of the main results of \cite{dGMS} concerns the construction
of various bases of the module $G_0(\log\,D)$ (hence, of all the other
modules given above), such that the connection takes a particularly simple form.
This can be summarized as follows.
\begin{proposition}[{\cite[proposition 4.5(iii)]{dGMS}}]\label{prop:recalldGMS}
There is a $\dC[\theta,t]$-basis $\underline{\omega}^{(1)}=(\omega_1,\ldots,\omega_n)$ of $G_0(\log\,D)$ such
that
\begin{equation}\label{eq:Basis}
\nabla(\underline{\omega}^{(1)}) = \underline{\omega}^{(1)} \cdot
\left[
(A_0\frac{1}{\theta}+A_\infty)\frac{d\theta}{\theta}+(-A_0\frac{1}{\theta}+A'_\infty)\frac{dt}{nt}
\right]
\end{equation}
where
$$
A_0:=
  \begin{pmatrix}
    0 & 0&\ldots &0& c\cdot t\\
    -1   &0&\ldots&0& 0\\
    \vdots & \vdots & \ddots & \vdots &\vdots \\
    0 &0 & \ldots &0& 0\\
    0&0& \ldots& -1&0
  \end{pmatrix},
$$
$A_\infty=\diag(\nu_1,\ldots,\nu_n)$ and
$A'_\infty:=\diag(0,1,\ldots,n-1)-A_\infty$.
The constant $c\in \dC$ is defined by the equation
$$
f^n=-c\cdot h + \sum_{i=1}^{n-1}\xi_i(f)\cdot k_i,
$$
where $\xi_i\in\Theta_{V/T}(\log\,D)$ and $k_i\in\cO_V$ are
homogenous polynomials of degree $n-1$.
Here
$$
\omega_1=\iota_E\frac{\vol}{h}=n\frac{\vol}{dh},
$$
where $\vol=dx_1\wedge\ldots\wedge dx_n$.
In particular, $(G(*D),\nabla)$ is flat.
Moreover, $(\nu_1,\ldots,\nu_n)$ is the spectrum at ($\theta=$) infinity
of the restriction of $G_0(\log\,D)$ to $t=0$.
\end{proposition}
The following obvious consequence will be used in the next section.
\begin{corollary}
Consider the inclusion $j:\{1\}\times T\hookrightarrow \dC\times T$ and
the restriction $(G_1(*D),\nabla):=j^*(G(*D),\nabla)$.
This is a meromorphic bundle on $T$ with connection (the only pole
being at $0\in T$), and can thus be seen
as a coherent and holonomic left $\Weyl$-module.
Then we have an isomorphism of left $\Weyl$-modules
$$
\varphi:\Weyl/(b_{G_1(\log\,D)}(t\partial_t)+\frac{c}{n^n}\cdot t) \stackrel{\cong}{\longrightarrow} (G_1(*D),\nabla),
$$
where $b_{G_1(\log\,D)}$ is the spectral polynomial of $G_1(\log\,D):=j^*(G(\log\,D),\nabla)$ at $0\in T$, i.e.,
$$
b_{G_1(\log\,D)}(s):=\prod_{i=1}^n \left(s-\frac{i-1-\nu_i}{n}\right),
$$
and where $\varphi(1)=\omega_1$. In particular, $\omega_1$ satisfies
the functional equation $b_{G_1(\log\,D)}(t\partial_t)\omega_1=-c/n^n\cdot t\omega_1$ in $G_1(*D)$.
\end{corollary}

We note the following easy consequence from the definitions, which was not stated in \cite{dGMS}.
\begin{lemma}\label{lemIndependenceLinearForm}
Let $f_1,f_2\in V^*\backslash D^*$ be two
generic linear forms.
Denote by $(G_i(\log\,D),\nabla_i)$ the family of Brieskorn lattices
attached to the pair $(f_i,h)$, $i=1,2$.
Then $\varphi_{c'}^*(G_1(\log\,D),\nabla_1) \cong (G_2(\log\,D),\nabla_2)$, where
$\varphi_{c'}:\dC\times T\rightarrow \dC\times T$ is defined as
$\varphi(\theta,t)=(\theta,c'\cdot t)$ for some $c'\in\dC^*$.
\end{lemma}
\begin{proof}
By definition, the complement of $D^*$ in $V^*$ is an (open) orbit
of the dual action of $G$, hence, there is $g\in G$ with $g(f_1)=f_2$.
Then $g(h)=\chi_h(g)\cdot h$ and it follows that
$\varphi_{\chi_h(g)}^*(G_1(\log\,D),\nabla_1) \cong (G_2(\log\,D),\nabla_2)$.
\end{proof}


In order to relate the above objects to the Bernstein polynomial of $h$, we recall how the Gau\ss-Manin system
$G(*D)$, seen as a left $\dC[\theta,t]\langle\partial_\theta,\partial_t\rangle$-module is obtained as a direct image
of a differential system on $V$. A similar reasoning as in the next lemma can be found in \cite[proposition 2.7]{DS}.
We consider $f$ as a morphism $f:V\rightarrow R=\Spec\dC[r]$, and put
$\Phi:=(f,h):V\rightarrow R\times T$.
\begin{lemma}\label{lem:CompDirectImage}
Let $\Phi_+ \cO_V(*D)$ be the (algebraic) direct image complex
of the holonomic $\cD_V$-module $\cO_V(*D)$. Then
\begin{enumerate}
\item
The cohomology sheaves of $\Phi_+\cO_V(*D)$ are $\cD_{R \times T}$-coherent, holonomic and regular.
\item
$\Phi_+\cO_V(*D)$, seen in $\cD^b(\cD_{R\times T /T})$ is represented by
$$
\left(\Phi_*\Omega^{\bullet+n-1}_{V/T}(*D)[\partial_r],d-(df\wedge-\otimes \partial_r)\right).
$$
Under the isomorphism
\begin{equation}\label{eq:GM-DirectImage}
\cH^0(\Phi_+\cO(*D)) \cong \frac{\Phi_*\Omega^{n-1}_{V/T}(*D)[\partial_r]}{(d-df\wedge-\otimes \partial_r)\Phi_*\Omega^{n-2}_{V/T}(*D)[\partial_r]}
\end{equation}
the action of $\partial_t$ on a class $[\omega]\in \cH^0(\Phi_+\cO_V(*D))$ represented
by a form $\omega\in \cH^{n-1}_{V/T}(*D)$ is given by
\begin{equation}
\label{eq:ConnectionParameter-BeforeFL}
\partial_t([\omega]) := \frac{1}{n\cdot t}\left([\Lie_E(\omega)]-[\Lie_E(f)\omega\otimes\partial_r]\right)
\end{equation}
\item
Put $M:=\dH^0(R\times T, \Phi_+\cO(*D)))$. Denote by $\widehat{M}$ the partial Fourier-Laplace transformation with respect to $r$
of $M$, i.e., $\widehat{M}=M$ as $\dC$-vector spaces, and we define an structure of a $\dC[\tau,t]\langle\partial_\tau,\partial_t\rangle$-module
on $\widehat{M}$ by $\tau\cdot:=\partial_r$ and $\partial_\tau:=-r\cdot$. Then $\widehat{M}$ is $\dC[\tau,t]\langle\partial_\tau,\partial_t\rangle$-holonomic,
with singularities at $\tau=\{0,\infty\}$ and $t=\{0,\infty\}$ at most, regular along $\{t=0\} \cup \{\tau=\infty\}$.
Moreover, by putting $\theta=\tau^{-1}$, the localized Fourier-Laplace transformation $\widehat{M}[\tau^{-1}]$  of $M$ is isomorphic to $G(*D)$
as a meromorphic vector bundle with connection (the one on $G(*D)$ being given by formula \eqref{eq:ConnectionOperator}).
\item
The restriction $(G_1(*D),\nabla)$ is isomorphic (as a left $\Weyl$-module) to $\dH^0(T, h_+\cO_V(*D)e^{-f})$, where
$\cO_V(*D)e^{-f}$ is the tensor product of $\cO_V(*D)$ with a rank one $\cO_V$-module formally generated
by $e^{-f}$, i.e., $\cO_V(*D)e^{-f}\cong \cO_V(*D)$ as $\cO_V$-modules, and the differential
of $\cO_V(*D)e^{-f}$ (i.e., the operator defining its $\cD_V$-module structure) is given
by $d_f:=d-df\wedge$.
\end{enumerate}
\end{lemma}
\begin{proof}
\begin{enumerate}
\item see \cite[theor\`eme 9.0-8.]{Mebkhout}
\item
It is well known that for any left $\cD_V$-module $\cM$, the direct image complex $\Phi_+ \cM$ is represented by
$$
\left(\bR\Phi_*\Omega^{n+\bullet}_V(\cM)[\partial_r,\partial_t],d-(df\wedge-\otimes\partial_r)-(dh\wedge-\otimes\partial_t)\right).
$$
Putting $\cM=\cO_V(*D)$ and using that $\Phi$ is affine, we consider the double complex
$$
E^{p,q}:=\left((\Phi_*\Omega^{p+q}_V(*D)[\partial_r])\otimes\partial_t^q,d-(df\wedge-\otimes\partial_r),-(dh\wedge-\otimes\partial_t )\right)
$$
whose total cohomology is $\cH^{p+q-n}(\Phi_+(\cO_V(*D)))$. The morphism $h$ is smooth restricted to $V\backslash D$, hence,
the Koszul complex $(\Omega^\bullet(*D),dh \wedge)$ is acyclic.
Therefore the second spectral sequence associated to the above double complex degenerates at the $E_2$-term and the isomorphism
$$
\Omega^{p-1}_{V/T}(*D) \stackrel{dh\wedge}{\longrightarrow} {\cK\!}er\left(\Omega^p_V(*D) \stackrel{dh\wedge}{\longrightarrow} \Omega^{p+1}_V(*D)\right)
$$
yields the above quasi-isomorphism.

In order to prove the formula for the action of $\partial_t$, notice that given a class $[\omega]$ defined
by a relative $n$-form $\omega\in\Omega^{n-1}_{V/T}(*D)$, the class corresponding to it in
$\cH^0(\Phi_+\cO_V(*D))$ is $[dh\wedge\omega]$. By definition, we have
$\partial_t([dh\wedge\omega])=[dh\wedge\omega\otimes\partial_t]$. This class is
equal in $\cH^0(\Phi_+\cO_V(*D))$ to $[d\omega-df\wedge\omega\otimes\partial_r]$.
It follows that under the isomorphism \eqref{eq:GM-DirectImage}, this equals
$[\frac{d\omega}{dh}-\frac{df\wedge\omega}{dh}\otimes\partial_r]\in
\Phi_*\Omega^{n-1}_{V/T}(*D)[\partial_r]/(d-df\wedge-\otimes\partial_r)\Phi_*\Omega^{n-2}_{V/T}(*D)$.
Now notice that as $h$ is smooth outside $D$, there is a vector field
$X\in\Theta_V(*D)$ which lifts $\partial_t\in\Theta_T$. Then we have that
$d\omega/dh=\iota_Xd\omega$ and $(df\wedge\omega)/dh=\iota_X(df)\wedge\omega-df\wedge\iota_X\omega$.
Putting this together and using once again the relation in
the quotient $\Phi_*\Omega^{n-1}_{V/T}(*D)[\partial_r]/(d-df\wedge-\otimes\partial_r)\Phi_*\Omega^{n-2}_{V/T}(*D)$
one arrives at the formula
$$
\partial_t[\omega]=[\Lie_X\omega]-[\Lie_X(f)\omega]\otimes\partial_r.
$$
Now the result follows as the meromorphic vector field $X$ can be taken to be
$E/(n\cdot h)$, due to the homogeneity of $h$.
\item
This is obvious from the last point: Fourier-Laplace transformation and localization along $\tau=0$
transforms formula \eqref{eq:GM-DirectImage} into
the defining equation \eqref{eq:GM-Brieskorn}
of $G(*D)$ and formula \eqref{eq:ConnectionParameter-BeforeFL} obviously
corresponds to the second part of formula \eqref{eq:ConnectionOperator}. The statements about regularity follows form the
general considerations in \cite[theorem 1.11]{DS}.
\item
By definition, $h_+\cO_V(*D)e^{-f}$ is represented by the complex
$$
(h_*\Omega^{n+\bullet}_V(*D)[\partial_t],d-df\wedge-(dh\wedge-\otimes\partial_t)).
$$
The same argument as above shows that this is quasi-isomorphic to
$$
(h_*\Omega^{n-1+\bullet}_{V/T}(*D),d-df\wedge).
$$
Now it is clear that
$$
\dH^0(T,h_*\Omega^{n-1+\bullet}_{V/T}(*D),d-df\wedge) = \frac{H^0(\Omega^{n-1}_{V/T}(*D))}{(d-df\wedge)H^0(V,\Omega^{n-2}_{V/T}(*D)))} = G_1(*D)
$$
\end{enumerate}
\end{proof}

\textbf{Remark:} Direct images of regular holonomic modules by a morphism consisting of two polynomials
occur in the work of C.~Roucairol (see \cite{Rouc1}, \cite{Rouc2} and \cite{Rouc3}). She also studied
direct images of twisted modules, i.e., $h_+(\cM e^{-f})$. However, we will not need her results directly
as the computations from \cite{dGMS} (i.e. proposition
\ref{prop:recalldGMS} above) give already very precise information about
these direct images for a
pair $(f,h)$, with $h$ reductive linear free and $f$ linear and generic.

\section{Bernstein Polynomials}
\label{sec:Bernstein}

We first give the definition of the Bernstein polynomial through
the classical functional equation. Next we recall how this can be rephrased
using the general theory of $V$-filtrations. This will be useful
in proving the main result. Finally, we state and prove the relation between
the roots of the Bernstein polynomial of a defining equation $h$ for a linear
free divisor and the residue eigenvalues of the family of Gau\ss-Manin-systems
introduced in section \ref{sec:LFDandGM}.

The following classical statement is due to Bernstein (see, \cite{Bernstein}).
\begin{theorem}\label{theo:BernsteinLFD}
Let $h\in \cO_V$ be any function, then there is a polynomial $B\in\dC[s]$ and a differential
operator $P(x_i,\partial_{x_i},s)\in\cD_V[s]$ such that
$$
P(x_i,\partial_{x_i},s) h^{s+1} = B(s) h^s
$$
All polynomials $B(s)\in\dC[s]$ having this property form an ideal in $\dC[s]$,
and we denote by $b_h(s)$ the unitary generator of this ideal. $b_h(s)$ is called
the Bernstein polynomial of $h$.
\end{theorem}
If $h$ defines a linear free divisor, then the theory of pre-homogenous vector spaces
shows that the functional equation defining $b_h(s)$ is of a particular type.
\begin{theorem}[\cite{sk}, \cite{Gyoja1}, \cite{GS}]\label{theo:BernsteinPolLinFree}
Let $D=h^{-1}(0)$ be a reductive linear free divisor, then the operator $P$ appearing in
Bernstein's functional equation is given by $P:=h^*(\partial_{x_1},\ldots,\partial_{x_1})$
(remember that $h^*(\underline{y})=\overline{h(\overline{\underline{y}})}$, where $x_i$ are the unitary coordinates
and $y_i$ are their duals).
In particular, it is an element of $\dC\langle\partial_{x_1},\ldots,\partial_{x_n}\rangle$.
Moreover, the degree of $b_h(s)$ is equal to $n$ and the roots of $b_h(s)$ are contained
in the open interval $(-2,0)$ and are symmetric around $-1$. In particular, $-1$ is the
only integer root.
\end{theorem}

The following classical reformulation of the definition of the Bernstein polynomial will
be useful in the sequel.

Consider the ring $\dC[s,h^{-1}]$, and denote by
$M[h^{-1}]:=\dC[s,h^{-1}] h^s$ the rank one $\dC[s,h^{-1}]$-module generated by the symbol $h^s$.
Define an action of $\cD_V$, the ring of algebraic differential operators on $V$ on
$M[h^{-1}]$ by putting $\partial_{x_i}(g\cdot h^s):=\partial_{x_i}(g)\cdot h^s+g\cdot s \cdot h^{-1}\partial_{x_i}(h)h^s$.
This action extends naturally to an $\cD_V[s]$-action. Let $M$ the $\cD_V[s]$-submodule of $M[h^{-1}]$ generated
by $h^s$. Define an action of $t$ on $M[h^{-1}]$ by putting $t(g(s)\cdot h^s):=g(s+1)\cdot h\cdot h^{s+1}$.
Then $b_h(s)$ is the minimal polynomial of the action of $s$ on the quotient $M/tM$.

This definition can be rephrased once more using the theory of $V$-filtrations on $\cD$-modules. Without reviewing the details
of the theory, we recall the following facts (see, e.g. \cite[section 4]{MebkhMais})
\begin{deflemma}
Let $X$ be any smooth algebraic variety, and $Y\subset X$ a smooth hypersurface defined by an ideal sheaf
$I\subset\cO_X$. We denote by $t\in\cO_X$ a local generator of $\cI$.
\begin{enumerate}
\item
Let $\cD_X$ be the sheaf of algebraic differential operators, then define
$$
\begin{array}{rcl}
V_k\cD_X:=\left\{ P\in \cD_X\,|\, P(I^j)\subset I^{j-k}\right\}
\end{array}
$$
For any left $\cD_X$-module $\cM$, a V-filtration on $\cM$ is an increasing filtration
$U_\bullet\cM$ compatible with $V_\bullet\cD_X$.
\item
A V-filtration $U_\bullet\cM$ on a left $\cD_X$-module $\cM$ is good iff the \emph{Rees}-module
$\oplus z^k U_k\cM$ is $\cR_V\cD_X$-coherent, where $\cR_V\cD_X:=\oplus_k z^k V_k\cD_X$.
\item
A good V-filtration $U_\bullet\cM$ is said to have a Bernstein polynomial iff there is
a non-zero polynomial $b(s)\in\dC[s]$ such that for all $k\in\dZ$, we have
$b(-\partial_t t + k)U_k\cM\subset U_{k-1}\cM$.
\item
A coherent $\cD_X$-module $\cM$ is called specializable iff locally there exists a good $V$-filtration
$U_\bullet\cM$ having a Bernstein polynomial. Equivalently, for any local section $m\in \cM$
there is a non-zero polynomial $b_m(s)$ (the Bernstein polynomial of $m$) such that $b_m(-\partial_t t)\,m\in V_{-1}\cD_X\cdot m$.
\item
A holonomic $\cD_X$-module is specializable along any smooth hypersurface $Y$.
\end{enumerate}
\end{deflemma}
The following evident corollary gives an example of a $V$-filtration that
will be used later.
\begin{corollary}
Consider the left $\Weyl$-module $G_1(*D)$ from above.
Then
$$
U_k G_1(*D):=V_k \Weyl \cdot G_1(\log\,D)
$$
defines
a good V-filtration on
$G_1(*D)$, whose Bernstein polynomial is exactly $b_{G_1(\log\,D)}(s)$.
Moreover, we have
$$
U_0 G_1(*D) = G_1(\log\,D) = V_0 \Weyl \cdot \omega_1.
$$
\end{corollary}

We will also use $V$-filtrations for $\cD_{T\times V}$-modules.
The following result is well known, see, e.g., \cite{Mal6}, \cite[lemme 4.4-1]{MebkhMais}.
\begin{lemma}\label{lem:BernsteinVFilt}
\begin{enumerate}
\item
Let $h\in \cO_V$ an arbitrary function, seen as a morphism $h:V\rightarrow T$.
Denote by $i_h:V\hookrightarrow T\times V$ the graph embedding, with image $\Gamma_h$. Put
$\cN:=(i_h)_+\cO_V$, then $\cN \cong \cO_V[\partial_t] \cong \cO_{T\times V}(*\Gamma_h)/\cO_{T\times V}\cong \cD_{T\times V}\delta(t-h)$.
A good $V$-filtration with respect to the hypersurface $\{0\}\times V$ on $\cN$ is defined by putting, for all $k \in \dZ$, $U_k \cN:=V_k\cD_{T\times V}\delta(t-h)$.
This $V$-filtration admits a Bernstein
polynomial (namely, a Bernstein polynomial for the section $\delta(t-h)$), which is exactly the polynomial $b_h(s)$. We denote, as in \cite{Mal6}, by
$\cM$ the $V_0\cD_{T\times V}$-module $U_0 \cN$.
\item
The direct image $(i_h)_+\cO_V(*D)$ is the localization of both $\cN$ and $\cM$ along $t=0$, and is thus
denoted by $\cM[t^{-1}]$. As $\cN$ has no $t$-torsion, we have an exact sequence
$$
0\longrightarrow \cN \longrightarrow \cM[t^{-1}] \longrightarrow C \longrightarrow 0.
$$
where $C$ is a $\cD_{T\times V}$-module.
A Bernstein
polynomial for a local section $m\in \cN$ is also a Bernstein polynomial for
$m$, seen as a local section in $\cM[t^{-1}]$.

\end{enumerate}
\end{lemma}

We can now state and prove the main result of this paper.
\begin{theorem}\label{theo:BernsteinEqualResidueGM}
Let $D=h^{-1}(0)$ be reductive linear free divisor and $f\in V^*$ be generic.
Consider the family of Gau\ss-Manin systems $G(*D)$, the
logarithmic extension $G(\log\,D)$
and the restrictions $G_1(\log\,D)\subset G_1(*D)$ from above.
Then we have that $b_h(s)=b_{G_1(\log\,D)}(s+1)$ (recall that $b_{G_1(\log\,D)}(s)$ is the spectral
polynomial of $G_1(\log\,D))$.
\end{theorem}
In order to prove this result, we start with a preliminary lemma.
\begin{lemma}
Let $\cM[t^{-1}]:=(i_h)_+\cO_V(*D)$ as above. Consider the
twisted module $(i_h)_+\cO_V(*D)e^{-f}$. Then the
section $\delta(t-h)e^{-f}\in(i_h)_+\cO_V(*D)e^{-f}$
admits $b_h(s)$ as a Bernstein polynomial, with associated
functional equation
\begin{equation}\label{eq:twistedBernstein}
\left(t\cdot h^*(\partial_{x_i}+a_i)-b_h(-\partial_t t)\right)\delta(t-h)e^{-f}=0,
\end{equation}
where $f=\sum_{i=1}^n a_ix_i$.
\end{lemma}
\begin{proof}
By lemma \ref{lem:BernsteinVFilt}, $b_h(s)$ is the minimal
polynomial of $-\partial_t t$ on
$
\frac{\cD_V[t\partial_t]\delta(t-h)}{t\cD_V[t\partial_t]\delta(t-h)}
$.
In particular, by theorem \ref{theo:BernsteinPolLinFree}, the functional equation
$$
\left(t\cdot h^*(\partial_{x_i})-b_h(-\partial_t t)\right)\delta(t-h)=0
$$
holds in $(i_h)_+\cO_V(*D)$. Then it follows directly from the definition of the
twisted module $\cO_V(*D)e^{-f}$ that the
functional equation \eqref{eq:twistedBernstein} from above holds in $(i_h)_+\cO_V(*D)e^{-f}$.
Now suppose that there
is another equation
$$
\left(t\cdot \widetilde{P}(x_i,\partial_{x_i},-\partial_t t)-\widetilde{B}(-\partial_t t)\right)\delta(t-h)e^{-f}=0,
$$
where $\widetilde{P}\in\cD_V[s]$ and $\widetilde{B}(s)\in\dC[s]$ with $\deg(\widetilde{B})<\deg(b_h)$.
Then we obtain the equation
$$
\left(t\cdot \widetilde{P}(x_i,\partial_{x_i}-a_i,-\partial_t t)-\widetilde{B}(-\partial_t t)\right)\delta(t-h)=0
$$
in $(i_h)_+\cO_V(*D)$, which contradicts the minimality of $b_h(s)$.
\end{proof}

\begin{proof}[Proof of the theorem]

We consider, as in the last lemma, the $\cD_{T\times V}$-module
$$
(i_h)_+\cO_V(*D)e^{-f}\cong (i_h)_*\cO_V(*D)e^{-f}[\partial_t]
$$
and the $\cD_V[t\partial_t]$-submodule
generated (over $\cD_V[t\partial_t]$) by $\delta(t-h)e^{-f}$. The direct image
$h_+(\cO_V(*D)e^{-f})$ is obtained in the standard way from
 $(i_h)_+\cO_V(*D)e^{-f}$ as the relative de Rham complex of the
projection $p_1:T\times V\rightarrow T$. In other words, we have
$$
\cH^i(h_+(\cO_V(*D)e^{-f}))  =  \cH^i((p_1)_*\DR^{n+\bullet}_{T\times V/T}((i_h)_*\cO_V(*D)e^{-f}[\partial_t]))
$$
Considering $(i_h)_+\cO_V(*D)e^{-f}$ as a $\cD_V$-module only, we
thus have
$$
\cH^i(h_+(\cO_V(*D)e^{-f}))=h_*\cH^i(\DR^{n+\bullet}_V\cO_V(*D)e^{-f}[\partial_t]).
$$
Now it is well known (see, e.g., \cite[proposition 2.1]{Mal6} or \cite[proposition 2.2.10]{Bj2}),
that for any left $\cD_V$-module $\cL$, the de Rham complex
$\DR_V^\bullet(\cL)$ represents the (shifted) derived tensor product
$\Omega^n_V\stackrel{\dL}{\otimes}_{\cD_V}\cL[-n]$, in particular, we have
$$
\cH^n(DR^\bullet_V(\cL)) \cong \Omega^n_V\otimes_{\cD_V}\cL.
$$
It follows that
\begin{equation}\label{eq:TopDeRham}
\cH^0(h_+(\cO_V(*D)e^{-f})) \cong h_*\left(\Omega^n_V\otimes_{\cD_V}(\cO_V(*D)e^{-f}[\partial_t])\right).
\end{equation}
so that, taking global sections and considering again the isomorphism
from lemma \ref{lem:CompDirectImage}, 4., we obtain
$$
H^0(V, \Omega^n_V\otimes_{\cD_V}(\cO_V(*D)e^{-f}[\partial_t])) \cong G_1(*D)
$$
Notice that the section $\vol\otimes\delta(t-h)e^{-f}$ is mapped to
the section $\omega_1/n=\vol/dh$ under this isomorphism.

From the equation
$\left(t\cdot h^*(\partial_{x_i}+a_i)-b_h(-\partial_t t)\right)\delta(t-h)e^{-f}=0$
in $\cO_V(*D)e^{-f}[\partial_t]$ (equation \eqref{eq:twistedBernstein}) we deduce that the element
$\vol\otimes \left(t\cdot h^*(\partial_{x_i}+a_i)-b_h(-\partial_t t)\right)\delta(t-h)e^{-f}$
is zero in $h_*(\Omega^n_V\otimes_{\cD_V}(\cO_V(*D)e^{-f}[\partial_t]))$.
%
Hence
$$
t\cdot (h^*(\partial_{x_i}+a_i)(\vol) )\otimes \delta(t-h)e^{-f} = b_h(-\partial_t t)(\vol\otimes \delta(t-h)e^{-f})
$$
holds in $h_*\left(\Omega^n_V\otimes_{\cD_V}(\cO_V(*D)e^{-f}[\partial_t])\right)$, where the operator $h^*(\partial_{x_i}+a_i)$ acts on $\vol$
by the right $\cD_V$-action on $\Omega^n_V$. Now develop the polynomial $h^*(y_i+a_i)$ as
$h^*(y_i+a_i) = \sum_{1 \leq |I|\leq n} a_I y^I + h^*(a_i)$, then
$$
h^*(\partial_{x_i}+a_i) = \sum_{1\leq|I|\leq n} a_I \partial_{x_1}^{i_1} \ldots \partial_{x_n}^{i_n} + h^*(a_i)
$$
and the action $h^*(\partial_{x_i}+a_i)(\vol)$ is given by
$$
\sum_{1\leq|I|\leq n} a_I \left(\underbrace{\Lie_{\partial_{x_1}}\cdots\Lie_{\partial_{x_1}}}_{i_1} \cdots
\underbrace{\Lie_{\partial_{x_n}}\cdots\Lie_{\partial_{x_n}}}_{i_n} \right)(\vol)+ h^*(a_i)\cdot\vol
$$
But obviously $\Lie_{\partial_{x_i}}\vol=0$ for any $i\in\{1,\ldots,n\}$, so that finally we see that
the section $\vol\otimes\delta(t-h)e^{-f}$ of $h_*\left(\Omega^n_V\otimes_{\cD_V}(\cO_V(*D)e^{-f}[\partial_t])\right)$ is annihilated
by $h^*(a_i)\cdot t-b_h(-\partial_t t)$. It follows that $b_h(-\partial_t t)$ sends $U_0 G(*D) = V_0\Weyl \omega_1$
into $U_{-1}G(*D)$, hence, we have $b_{G_1(\log\,D)}(s+1)|b_h(s)$.
Now the theorem
follows as both $b_h$ and $b_{G_1(\log\,D)}$ are of degree $n$.
\end{proof}

\section{Consequences and Examples}
\label{sec:ConsqExamples}

\begin{definition}
Let $D$ be a reductive linear free divisor with defining equation $h\in\cO_V$ and
$f\in V^*$ a generic linear form. Consider, as in the last section, the
logarithmic extension $G_0(\log\,D)$ of the family of Brieskorn lattices $G_0(*D)$
attached to $(f,h)$.
We define the \emph{logarithmic Brieskorn lattice} of $h$ to be the restriction $G_0(h):=i^*(G_0(\log\,D),\nabla)$,
where $i:\dC\times\{0\}\hookrightarrow \dC\times T$.
\end{definition}
Notice that it follows from lemma \ref{lemIndependenceLinearForm} that
$G_0(h)$ is independent of the choice of $f$ in $V^*\backslash D^*$, so that
it makes sense to speak about \emph{the} logarithmic Brieskorn lattice of $h$.


The next result, which is an easy consequence of theorem \ref{theo:BernsteinEqualResidueGM}, can be considered as a variant
of the corresponding classical statement of Malgrange (\cite{Mal6}) for the isolated singularity case.
\begin{theorem}\label{theo:AdaptMalgrange}
Let $(G_0(h),\nabla)$ be the logarithmic Brieskorn lattice of a reductive linear free divisor $D$. Then
$\nabla$ is regular singular at $\theta=0$. Consider the saturation
$\widetilde{G}_0(h):=\sum_{k\geq 0} (\nabla_{\theta\partial_\theta})^kG_0(h)$, which
has a logarithmic pole at $\theta=0$.  Let $b_{\widetilde{G}_0(h)}(s)$ be the minimal polynomial of the residue endomorphism of $\nabla_\theta$ on $\widetilde{G}_0(h)$. Then $b_{\widetilde{G}_0(h)}(n(s+1)) = b_h(s)$.

\end{theorem}
\begin{proof}
The regularity follows easily from the particular form of the connection matrix \eqref{eq:Basis}. Namely,
$G_0(h)$ is the Fourier-Laplace transformation of a regular $\dC[r]\langle \partial_r\rangle$-module, hence, its
regularity is equivalent to the nilpotency of the polar part of the connection matrix, which is obviously
the case here, by putting $t=0$ in $A_0$. Now the saturation of $G_0(h)$ is easy to calculate:
We put $\widetilde{\omega}_i:=\theta^{1-i}\omega_i$, then $G(\log\,D)=
\oplus_{i=1}^n\dC[\theta,\theta^{-1},t]\widetilde{\omega}_i$,
but $G_0(\log\,D)\subsetneq\oplus_{i=1}^n\dC[\theta,t]\widetilde{\omega}_i$.
It is evident
that $\widetilde{G}_0(h)=\oplus_{i=1}^n\dC[\theta]\widetilde{\omega}_i$, in particular,
this module is invariant under $\theta\nabla_\theta$, i.e., logarithmic at $\theta=0$.
We have $(\theta\partial_\theta)\underline{\widetilde{\omega}}=\underline{\widetilde{\omega}}
\cdot(\widetilde{A}_0+\diag(\{1-i+\nu_i\}_{i=1,\ldots,n})$, where $\widetilde{A}_0:=(A_0)_{|t=0}$.
We see by theorem \ref{theo:BernsteinEqualResidueGM} that the
residue eigenvalues of $\nabla_\theta$ at $\theta=0$ are the roots of the Bernstein polynomial
of $h$ after dividing by $n$ and shift by $-1$, and moreover that the
residue endomorphism is regular (i.e., its minimal and characteristic polynomial
coincide), as it has a cyclic generator. This proves the theorem.
\end{proof}
\textbf{Remark: } One may ask what the meaning of the rescaling by $n$
occurring in $b_{\widetilde{G}_0(h)}(n(s+1))$ is.
The same kind of twist occurs in \cite[proposition 4.5i(v)]{dGMS}, where it is performed
on the base, i.e., where the pull-back $u^*(G(*D))$ with $u:\dC^2\rightarrow\dC\times T$, $(\theta,t')\mapsto(\theta,(t')^n)$
is considered, and where it is shown that after this pull-back, the resulting bundle has the
``rescaling property'', i.e., that it is invariant under $\nabla_{\theta\partial_\theta - t'\partial_{t'}}$.

\vspace*{0.5cm}

The following easy consequence is a somewhat reverse argumentation compared to
Malgrange's result, where the rationality of the roots of the
Bernstein polynomial was deduced from the (known) quasi-unipotency of the monodromy
acting on the cohomology of the Milnor fibre of an isolated hypersurface singularity. In our case, the rationality of the
roots of $b_h(s)$ is known, but we deduce information on the
(a priori unknown) monodromy of the logarithmic Brieskorn lattice $G_0(h)$. Moreover,
we can use the results of \cite{GS} to obtain a symmetry property of the spectrum
at infinity of the logarithmic Brieskorn lattice, which was conjectured in \cite[corollary 5.6]{dGMS}.
\begin{corollary}\label{cor:MondromySymmetryAtZero}
The monodromy of the logarithmic Brieskorn lattice, i.e. of the
local system associated to $G_0(h)[\theta^{-1}]:=G_0(h)\otimes_{\dC[\theta]}\dC[\theta,\theta^{-1}]$ is quasi-unipotent.
Moreover, let $\alpha_1,\ldots,\alpha_n$ be the spectral numbers of $G_0(h)$
at infinity (i.e., the numbers $\nu_i$ from proposition \ref{prop:recalldGMS}), written as a non-decreasing sequence. Then
$\alpha_i+\alpha_{n+1-i}=n-1$.
\end{corollary}
\begin{proof}
The eigenvalues of this monodromy are simply the exponentials of either
the numbers $\nu_i$ or $\nu_i':=i-1-\nu_i$ from proposition \ref{prop:recalldGMS} (or any other
integer shift of them). The numbers $\nu_i'$ are the roots of the
Bernstein polynomial of $h$ shifted by one, as shown in theorem \ref{theo:BernsteinEqualResidueGM}.
These are known to be rational by \cite{Ka1}. Similarly, if we denote
the roots of $b_h$ by $\alpha'_1,\ldots,\alpha'_n$, with $\alpha'_i\leq \alpha'_j$ if $i\leq j$,
then we know from
\cite[theorem 2.5.]{GS} that $\alpha'_i+\alpha'_{n+1-i}=-2$. From theorem
\ref{theo:BernsteinEqualResidueGM} and proposition \ref{prop:recalldGMS} we deduce
that $\alpha_j=(j-1)-\alpha'_j-1$ for any $j\in\{1,\ldots,n\}$, hence, $\alpha_i+\alpha_{n+1-i}=((i-1)-\alpha'_i-1)+((n+1-i)-1-\alpha'_{n+1-i}-1)=n-1$.
\end{proof}

We outline another consequence of the theorem \ref{theo:BernsteinEqualResidueGM}. Its interest is motivated
by comparing the situation considered here with the one where $f$ is still a generic linear form, but
$h$ is supposed to by an arbitrary monomial $h=\prod x_i^{w_i}$, i.e., non-reduced. The corresponding Gau\ss-Manin-systems
resp. Brieskorn lattices have been studied in \cite{DS2}, \cite{Dou3} and \cite{DM}. It is known that they are closely
related to the \emph{Mirror symmetry} phenomenon, i.e., one constructs a Frobenius structure on the
semi-universal unfolding of $f_{|h^{-1}(t)}$, $t\neq 0$ which is known to be isomorphic to
the orbifold quantum cohomology of the weighted projective spaces. For a linear free divisor
$D$, a similar construction of a Frobenius manifold has been carried out in \cite{dGMS}. Although
these are not a priori mirrors of some variety or orbifold, the following corollary shows an interesting
similarity with the case $h=\prod x_i^{w_i}$.

\begin{corollary}
The spectrum at $\theta=\infty$ of both $(G_0(h),\nabla)$ and $(G(*D),\nabla)$ contains a (non-trivial) block
of integer numbers $k,k+1,\ldots,n-1-k$ for some $k\in\{0,\ldots,n-1\}$.
\end{corollary}
\begin{proof}
For the spectrum of $(G_0(h),\nabla)$, this is obvious as this block corresponds to the root $-1$ of the Bernstein
polynomial $b_h(s)$. For the spectrum of $(G(*D),\nabla)$, one shows the same statement by analyzing the construction of
a good basis of $G(*D)$ from a good basis of $G_0(h)$ using algorithm 2 of \cite[lemma 4.11]{dGMS}.
\end{proof}
Notice that for the normal crossing case, the integer $k$ from above is equal to zero, i.e., the block
mentioned above is the whole spectrum. This is not true in general, hence, the Frobenius structures
constructed in \cite{dGMS} are not, a priori, mirrors of quantum cohomology algebras of orbifolds,
as zero is not, in general, an element of the spectrum. Still the analogy with the orbifold quantum
cohomology, i.e., the fact that there is a block of increasing integer spectral numbers corresponding
to the ``untwisted sector'' (see, e.g., \cite[section 2.1.]{Ir2}) is rather intriguing.

\vspace*{0.5cm}

\textbf{Examples of Bernstein polynomials:} We use the main result and the computations
of spectral numbers in \cite{dGMS} to obtain the roots of the Bernstein polynomials
for the following reductive linear free divisors. The definitions of the two last discriminants
can be found in \cite {GMNS}, example in 1.4(2) (this one is also called ``bracelet'') and \cite{sk}, proposition 11, respectively.\\

\begin{table}[h]
$$
\begin{array}{c|c}
\textup{\textbf{linear free divisor}} &  \textup{\textbf{Bernstein polynomial of \textit{h}}}\\ \hline \\
A_n\textup{ - quiver} & (s+1)^n \\ \hline \\
D_m\textup{ - quiver} & \left(s+\frac43\right)^{m-3}\cdot(s+1)^{2m-4}\cdot\left(s+\frac23\right)^{m-3} \\ \hline \\
E_6\textup{ - quiver} & (s+\frac75)\cdot(s+\frac43)^4\cdot(s+\frac65)\cdot(s+1)^{10}\cdot(s+\frac45)\cdot(s+\frac23)^4\cdot(s+\frac35) \\ \hline \\
\star_m\textup{ - quiver} &  \prod_{l=0}^{m-3}\left(s+\frac{2(m-1)-l}{m}\right)^{l+1}\cdot(s+1)^{2(m-1)}\cdot\prod_{l=0}^{m-3}\left(s+\frac{m-1-l}{m}\right)^{m-l-2}\\ \hline \\
\textup{discriminant in }S^3((\dC^2)^*) & \left(s+\frac76\right)\cdot\left(s+1\right)^2\cdot\left(s+\frac56\right) \\ \hline \\
\textup{discriminant of }\\\textup{Sl}(3,\dC)\times \Gl(2,\dC)\textup{ action}\\\textup{on Sym}(3,\dC)\times\textup{Sym}(3,\dC)& \left(s+\frac54\right)^2\cdot\left(s+\frac76\right)^2\cdot\left(s+1\right)^4\cdot
\left(s+\frac56\right)^2\cdot\left(s+\frac34\right)^2
\end{array}
$$
\caption{Bernstein polynomials for some examples of linear free divisors}
\label{tab:ResultsBernsteinPol}
\end{table}
Notice that the examples $E_6$ and the last two discriminants are obtained by direct
calculations in Singular (\cite{Singular}). On the other hand, the closed formulas for
the star quiver and the $D$-series follows from rather involved combinatorial
arguments, the details of which will appear in \cite{CompQuiver}.
The Bernstein polynomials for $D_4$ (which is equal to $\star_3$) and the
bracelet are also calculated in \cite{GS}. The one for $A_n$ is of course completely obvious
and well known. It would be of interest to complete these calculations
by the Bernstein polynomials of quiver representations for the highest roots of the Dynkin quivers $E_7$ and $E_8$,
however, this seems to be out of reach of computer algebra for the moment (remember from \cite{BM} that
the linear free divisors associated to these roots for $E_7$ resp. $E_8$ are of degree $46$ resp. $118$).

\vspace*{0.5cm}

Let us finish this note with a remark and a conjecture exploiting further the analogy with the case of an isolated hypersurface
singularity. We have seen that the theorem of Malgrange can be adapted for reductive linear free divisors using the logarithmic Brieskorn lattice from above.
The regularity of $(G_0(h),\nabla)$ at $\theta=0$
suggest to study the spectrum in the classical sense of Varchenko (i.e., at $\theta=0$) of this lattice.
We recall the definition and calculate two examples, in order to show that this spectrum contains additional information
not present in roots of the Bernstein polynomial, similarly to the case of isolated singularities.
\begin{definition}
Let $(\mathbb{E},\nabla)$ be a vector bundle on $\dC=\Spec\dC[\theta]$ equipped with a connection with a pole at zero of order two at most, which is regular singular.
The localization $\mathbb{M}:=E\otimes_{\dC[\theta]}\dC[\theta,\theta^{-1}]$ has the
structure of a holonomic $\dC[\theta]\langle\partial_\theta\rangle$-module with a regular singularity at $\theta=0$.
We suppose that the monodromy of its de Rham complex is quasi-unipotent.
Denote by $V^\bullet \mathbb{M}$ the canonical V-filtration on
$\mathbb{M}$ at $\theta=0$, indexed by $\dQ$. As this is a filtration by free $\dC[\theta]$-modules
(and not by free $\dC[\theta^{-1}]$-modules as the $V$-filtration at $\theta=\infty$), we
write it as a decreasing filtration. Define the spectrum of $(E,\nabla)$ to be
$$
\Sp(E,\nabla):=\sum_{\alpha \in \dQ} \frac{V^\alpha \mathbb{M} \cap \mathbb{E}}{V^\alpha \mathbb{M} \cap \theta\mathbb{E} + V^{>\alpha} \mathbb{M} \cap \mathbb{E}} \alpha \in \dZ[\dQ]
$$
where $V^{>\alpha} \mathbb{M} :=\cup_{\beta > \alpha} V^\alpha \mathbb{M}$.
\end{definition}

As an example, we consider the case of the normal crossing divisor
$D=\{h^{A_n}=\prod_{i=1}^n x_i=0\}$, which is the discriminant
in the representation space of the quiver $A_n$.
It was stated in \cite{dGMS} (but essentially well known before,
due to the relation of this example to the quantum cohomology
of the projective space $\dP^{n-1}$) that we have
$G_0(h^{A_n}):=\oplus_{i=1}^n \cO_{\dC\times\{0\}}\omega_i$, and
$$
\nabla(\underline{\omega}) = \underline{\omega} \cdot \left[\frac{\widetilde{A}_0}{\theta}+\diag(0,1,\ldots,n-1)\right]\frac{d\theta}{\theta},
$$
$\widetilde{A}_0:=(A_0)_{|t=0}$.
On the
other hand, we take up the
example of the star quiver with three exterior vertices studied in \cite[example 2.3(i)]{dGMS}.
Notice that this is exactly the quiver $D_4$. Here $D\subset V=\dC^6$, and $h^{\star_3}=h^{\star_3}_1\cdot h^{\star_3}_2 \cdot h^{\star_3}_3$, where
$$
h^{\star_3}_1=\left|\begin{array}{cc} a & b \\ d & e \end{array}\right|
\quad;\quad
h^{\star_3}_2=\left|\begin{array}{cc} a & c \\ d & f \end{array}\right|
\quad;\quad
h^{\star_3}_3=\left|\begin{array}{cc} b & c \\ e & f \end{array}\right|.
$$
Following the various algorithms of loc.cit used to obtain good basis, we have that
$G_0(h^{\star_3}):=\oplus_{i=1}^6 \cO_{\dC\times\{0\}}\omega_i$, and
$$
\nabla(\underline{\omega}) = \underline{\omega} \cdot \left[\frac{A_0}{\theta}+\diag(2,1,2,3,4,3)\right]\frac{d\theta}{\theta}
$$
Notice that this is the basis called $\underline{\omega}^{(2)}$ in loc.cit.
\begin{proposition}
\begin{enumerate}
\item
The spectrum at $\theta=0$ for $h^{A_n}$ is
$$
\Sp(G_0(h^{A_n}),\nabla)=(0,1,\ldots,n-1)\in\dZ[\dQ],
$$
hence, it is equal to the spectrum at $\theta=\infty$ of both
$(G_0(h),\nabla)$ and $(G(*D),\nabla)$ (so that in this case we do not
get more information from the spectrum at $\theta=0$ than those
contained in the roots of $b_h(s)$).
\item
The spectrum at $\theta=0$ for $h^{\star_3}$ is given by
$$
\Sp(G_0(h^{\star_3}),\nabla)=(-2,1,2,3,4,7)\in\dZ[\dQ],
$$
hence, different from $\Sp(G_0(h),\nabla)$ and not directly related
to $b_h(s)=(s+\frac43)(s+1)^4(s+\frac23)$.
\end{enumerate}
\end{proposition}
\begin{proof}
\begin{enumerate}
\item
One can calculate directly that $G_0(h^{A_n})$ can be generated by elementary sections,
which implies that $\Sp(G_0(h^{A_n}),\nabla)$ is equal to the spectrum at $\theta=\infty$,
i.e., $\Sp(G_0(h^{A_n}),\nabla)=(0,1,\ldots,n-1)$. However, this can also be obtained in a more abstract
way: For any linear free divisor $D$, the analytic object corresponding to the restriction of $G_0(*D)$ to $\dC \times (T\backslash\{0\})$ is known
(after a finite ramification of order $n$) to be a Sabbah orbit of TERP-structures
(see the remark after the proof of theorem \ref{theo:AdaptMalgrange} and \cite[proposition 4.5 (v)]{dGMS}). In the $A_n$-case, it is easy to see that
the extension $G_0(\log\,D)$ is exactly the extension ${_0}\cE$ considered in
\cite[proof of theorem 7.3 and lemma 6.11]{HS1} and the logarithmic Brieskorn lattice $G_0(h)$ is
isomorphic to the limit $\cG_0$ considered in loc.cit, proof of theorem 7.3 and lemma 6.12. It was
shown in the proof of theorem 7.3 of loc.cit. that $\cG_0$ is generated by elementary sections.
\item
In the $\star_3$-case, one cannot apply the previous reasoning. Hence a direct calculation
is necessary. We explain parts of it, leaving the details to the reader.
From the connection matrix given above we see that
$(\theta\partial_\theta)\omega_6=3\omega_6$, and
$(\theta\partial_\theta)\omega_5=4\omega_5-\theta^{-1}\omega_6$.
We make the Ansatz
$$
\omega_5=\alpha\theta^{-1}\omega_6+s_4
$$
where $s_4$ is a section of $G_0(h)[\theta^{-1}]$ satisfying
$(\theta\partial_\theta)(s_4)=4\cdot s_4$. We obtain
$$
(\theta\partial_\theta)\omega_5 = 2\alpha\theta^{-1}\omega_6+4s_4
\stackrel{!}{=}
(4\alpha-1)\theta^{-1}\omega_6+4s_4
$$
from which we conclude that
$\omega_5=\frac12\theta^{-1}\omega_6+s_4$. Similarly, the equation
$(\theta\partial_\theta)\omega_4=3\omega_4-\theta^{-1}\omega_5$
is satisfied by putting
$$
\omega_4=\frac18 \beta_1 \cdot \theta^{-2}\omega_6+ s_3
$$
where $s_3\in G_0(h)[\theta^{-1}]$ is a section satisfying
$(\theta\partial_\theta) s_3 = 3 s_3 + \theta^{-1} s_4$.
Continuing this way we see that the elements of our basis
$\underline{\omega}$ can be written as finite sums of elementary
sections in the following way:
$$
\begin{array}{c}
\omega_1  =
\frac{1}{128}\theta^{-5}\omega_6+\frac{1}{16}\theta^{-4}s_4+ \frac18
\theta^{-3}s_3+\frac14 \theta^{-2}s_2 + \frac12 \theta^{-1} s_1
+\widetilde{s}_2
\\ \\
\begin{array}{rclcrclcrcl}
\omega_2 & = & \frac{1}{32}\theta^{-4}\omega_6+s_1 & ; &
\omega_3 & = & \frac{1}{16}\theta^{-3}\omega_6-s_2 & ; &
\omega_4 & = & \frac18\theta^{-2}\omega_6+s_3
\end{array}
\\ \\
\begin{array}{rclcrcl}
\omega_5 & = & \frac12\theta^{-1}\omega_6+s_4 & ; &
\omega_6 & = & \omega_6
\end{array}
\end{array}
$$
where $s_1,s_2,s_3,s_4,\widetilde{s}_2$ are sections of
$G_0(h)[\theta^{-1}]$ satisfying
$$
\begin{array}{rclcrcl}
(\theta\partial_\theta) s_1 & = &   s_1 + \theta^{-1} s_2 & ; &
(\theta\partial_\theta) s_2 & = & 2 s_2 + \theta^{-1} s_3 \\ \\
(\theta\partial_\theta) s_3 & = & 3 s_3 + \theta^{-1} s_4 & ; &
(\theta\partial_\theta) s_4 & = & 4 s_4 \\ \\
(\theta\partial_\theta) \widetilde{s}_2 & = & 2 \widetilde{s}_4
\end{array}
$$
Now it is easy to calculate an upper triangular base change
yielding a good basis and to show that  the spectrum is
$$
\Sp(G_0(h^{\star_3}),\nabla)=(-2,1,2,3,4,7)\in\dZ[\dQ],
$$
as required.
\end{enumerate}
\end{proof}
Based on the computations of these examples, we state the following conjecture, which is
related to corollary \ref{cor:MondromySymmetryAtZero} as well as to \cite[conjecture 5.5.]{dGMS}.
\begin{conjecture}
Let $h$ be the defining equation of a reductive linear free divisor $D\subset V=\dC^n$. Then the
spectrum of its logarithmic Brieskorn lattice $(G_0(h),\nabla)$ at $\theta=0$ is symmetric around $\frac{n-1}{2}$.
\end{conjecture}
\textbf{Remark: }
There are several questions one may ask about the spectrum at
$\theta=0$. First, it is surprising that negative numbers
(even smaller than $-1$) occur in this spectrum. One might want to
understand the possibly range for the spectrum, as well
as the difference to the roots of $b_h$, when multiplied by $n$.
This should be compared to the results in \cite{HeSt}
for isolated singularities, in particular,
lemma 3.4 of loc.cit.

\bibliographystyle{amsalpha}

\providecommand{\bysame}{\leavevmode\hbox to3em{\hrulefill}\thinspace}
\providecommand{\MR}{\relax\ifhmode\unskip\space\fi MR }
\providecommand{\MRhref}[2]{%
  \href{http://www.ams.org/mathscinet-getitem?mr=#1}{#2}
}
\providecommand{\href}[2]{#2}

\vspace*{1cm}

\nd
Lehrstuhl f\"ur Mathematik VI \\
Institut f\"ur Mathematik\\
Universit\"at Mannheim,
A 5, 6 \\
68131 Mannheim\\
Germany

\vspace*{1cm}

\nd
Christian.Sevenheck@math.uni-mannheim.de

\end{document}